\documentclass{ifacconf}
\usepackage{natbib}            % you should have natbib.sty
\usepackage{graphicx} 
\usepackage[utf8]{inputenc}
\usepackage[T1]{fontenc}
\usepackage[retainorgcmds]{IEEEtrantools}
\usepackage{amssymb,amsmath}
\usepackage{url}
\usepackage{siunitx}
\usepackage{xcolor}

\graphicspath{{Figures/}}

\newcommand{\abs}[1]{\left\lvert#1\right\rvert}
\newcommand{\bigabs}[1]{\big\lvert#1\big\rvert}
\newcommand{\Bigabs}[1]{\Big\lvert#1\Big\rvert}
\newcommand{\inner}[2]{\langle#1,#2\rangle}
\newcommand{\SO}{\mathrm{SO(3)}}
\newcommand{\Mab}{\Bigl(\frac{l_\alpha}{k_\alpha}a_{i\times}^2 +\frac{l_\beta}{k_\beta}b_{i\times}^2\Bigr)}

\begin{document}

\begin{frontmatter}

\title{A global observer for attitude and gyro biases from vector measurements}
\thanks[footnoteinfo]{This work was supported by the French Agence Nationale de la Recherche through the ANR ASTRID SCAR project ``Sensory Control of Aerial Robots'' (ANR-12-ASTR-0033).}

\author[First]{Philippe Martin}
\author[First]{Ioannis Sarras}

\address[First]{Centre Automatique et Systèmes, MINES ParisTech, PSL Research University, Paris, France\\
	(e-mail: \{philippe.martin,ioannis.sarras\}@mines-paristech.fr).}

\begin{keyword}                      % Five to ten keywords,
Attitude estimation; nonlinear observer; guidance; navigation systems.
\end{keyword}                             % keyword list or with the
                                          % help of the Automatica
                                          % keyword wizard

\begin{abstract}
	We consider the classical problem of estimating the attitude and gyro biases of a rigid body from vector measurements and a triaxial rate gyro. We propose a simple ``geometry-free'' nonlinear observer with guaranteed uniform global asymptotic convergence and local exponential convergence; the stability analysis, which relies on a strict Lyapunov function, is rather simple. The excellent behavior of the observer is illustrated through a detailed numerical simulation.
\end{abstract}

\end{frontmatter}

\section{Introduction}\label{sec:intro}

Estimating the attitude of a rigid body from vector measurements has been for decades a problem of interest, because of its importance for a variety of technological applications such as satellites or unmanned aerial vehicles.

The attitude of the body can be described by the rotation matrix $R\in\SO$ from inertial to body axes. The measurement vectors $u_1,\cdots,u_n\in\Rset^3$ correspond to the expression in body axes of known vectors $U_1,\cdots,U_n\in\Rset^3$ which are constant in inertial axes, i.e., $u_k(t)=R^T(t)U_k$. The goal is to reconstruct the attitude at time~$t$ using only the knowledge of the measurement vectors until~$t$. The problem would be very easy if the measurements were perfect: indeed, using for instance only the two vectors $u_1(t)$ and $u_2(t)$ and noticing $R^T(x\times y)=R^Tx\times R^Ty$ since $R$ is a rotation matrix, we readily find
\begin{IEEEeqnarray*}{rCl}
	R^T(t) &=& R^T(t)\cdot\begin{pmatrix}U_1& U_2& U_1\times U_2\end{pmatrix}\cdot\begin{pmatrix}U_1& U_2& U_1\times U_2\end{pmatrix}^{-1}\\
	&=& \begin{pmatrix}u_1(t)& u_2(t)& u_1(t)\times u_2(t)\end{pmatrix}\cdot\begin{pmatrix}U_1& U_2& U_1\times U_2\end{pmatrix}^{-1}.
\end{IEEEeqnarray*} 

But in real situations, the measurement vectors are always corrupted at least by noise. Moreover, the $U_k$'s may possibly be not strictly constant: for instance a triaxial magnetometer measures the (locally) constant Earth magnetic field, but is easily perturbed by ferromagnetic masses and electromagnetic perturbations; similarly, a triaxial accelerometer can be considered as measuring the direction of gravity provided it is not undergoing a substantial acceleration (see e.g. \cite{MartiS2010ICRA} for a detailed discussion of this assumption and its consequences in the framework of quadrotor UAVs). 
That is why, despite the additional cost, it may be interesting to use a triaxial rate gyro to supplement the possibly not so good vector measurements.

The literature on attitude estimation from vector measurements can be broadly divided into three categories: i) optimization-based methods; ii) stochastic filtering; iii) nonlinear observers. Details on the various approaches can be found e.g. in the surveys \cite{CrassMC2007JGCD,ZamanTM2015arXiv} and the references therein.  

The first category sets the problem as the minimization of a cost function, and is usually referred to as Wahba's problem. The attitude is algebraically recovered at time~$t$ using only the measurements at time~$t$. No filtering is performed, and possibly available velocity information from rate gyros is not exploited.

The second category mainly hinges on Kalman filtering and its variants. Despite their many qualities, the drawback of those designs is that convergence cannot in general be guaranteed except for mild trajectories. Moreover the tuning is not completely obvious, and the computational cost may be too high for small embedded processors.

The third, and more recent, approach proposes nonlinear observers with a large guaranteed domain of convergence and a rather simple tuning through a few constant gains.
These observers can be designed: a) directly on $\SO$ (or the unit quaternion space), see e.g.~\cite{MahonHP2008TAC,MartiS2010CEP,VascoSO2008IFAC,GripFJS2012TAC}; b) or more recently, on $\Rset^{3\times3}$, i.e., deliberately ``forgetting'' the underlying geometry~\cite{BatisSO2012CEP,BatisSO2014AUT,GripFJS2015AUT}. Probably the best-known design is the so-called nonlinear complementary filter of~\cite{MahonHP2008TAC}; as noticed in~\cite{MartiS2007CDC}, it is a special case of so-called invariant observers~\cite{BoMaRo_tac09}.

In this paper, we propose a new observer of attitude and gyro biases from gyro measurements and (at least) two measurement vectors.  It also ``forgets'' the geometry of~$\SO$, which yields a very simple structure with a straightforward proof of its uniform global asymptotic convergence (notice the observer of~\cite{MahonHP2008TAC} is only quasi-globally convergent, with moreover a much more involved proof). It is easy to tune and implement, and demonstrates excellent performance in simulation. This observer can be seen as a modification and simplification of the linear cascaded observer proposed in~\cite{BatisSO2012CEP}, along with a much simpler proof of convergence based on a strict Lyapunov function. Compared to the global exponential stability result of~\cite{BatisSO2012CEP}, we only prove uniform global asymptotic stability plus local exponential stability. However, in the case where our observer reduces to the observer of~\cite{BatisSO2012CEP} by using directly non-filtered measurements in certain terms, our Lyapunov function again serves as a strict Lyapunov function, which provides a simpler and more direct alternative to the non-trivial stability analysis of~\cite{BatisSO2012CEP}.

The paper runs as follows: the model used to design the observer is described in section~\ref{sec:Models}; the observer is presented in section~\ref{sec:Observer}, and its convergence is proved; finally, section~\ref{sec:Simulations} illustrates the excellent behavior of the observer on a detailed numerical simulation.

\section{The design model}\label{sec:Models}
We consider a moving rigid body subjected to the angular velocity~$\omega$ (in body axes). Its orientation (from inertial to body axes) matrix $R\in\SO$ is related to~$\omega$ by
\begin{IEEEeqnarray}{rCl}
	\dot R &=& R\omega_\times\label{eq:R},
\end{IEEEeqnarray} 
where the skew-symmetric matrix $\omega_\times$ is defined by $\omega_\times x:=\omega\times x$ whatever the vector~$x$.

The rigid body is equipped with a triaxal rate gyro measuring the angular velocity~$\omega$, and two additional triaxial sensors (for example accelerometers, magnetometers or sun sensors) providing the measurements of two vectors $\alpha$ and~$\beta$. These two vectors correspond to the expression in body axes of two known independent vectors $\alpha_i$ and $\beta_i$ which are constant in inertial axes. In other words,
\begin{IEEEeqnarray*}{rCl}
	\alpha &:=&  R^T\alpha_i\\
	\beta &:=&  R^T\beta_i.
\end{IEEEeqnarray*}
Since $\alpha_i$, $\beta_i$ are constant, we readily find 
\begin{IEEEeqnarray*}{rCl}
	\dot\alpha &=& \alpha\times\omega\\
	\dot\beta &=& \beta\times\omega.
\end{IEEEeqnarray*}
As any sensor, the rate gyro is biased, and rather provides the measurement
\begin{IEEEeqnarray*}{rCl}
	\omega_m &:=& \omega+b,
\end{IEEEeqnarray*}
where $b$ is a slowly-varying (for instance with temperature) unknown bias. The effect of this bias on attitude estimation may be important when the observer gains are small, hence it is worth determining its value. But being not exactly constant, it cannot be calibrated offline and must be estimated online together with the attitude.

Our objective is to design an estimation scheme that can reconstruct online the orientation matrix $R(t)$ and the bias $b(t)$, using i) the measurements of the gyro and of the two vector sensors; ii) the knowledge of the constant vectors $\alpha_i$ and~$\beta_i$. The model on which the design will be based therefore consists of the dynamics
\begin{IEEEeqnarray}{rCl}
	\dot\alpha &=& \alpha\times\omega\label{eq:alpha}\\
	\dot\beta &=& \beta\times\omega\\
	\dot b &=& 0,
\end{IEEEeqnarray}
together with the measurements
\begin{IEEEeqnarray}{rCl}
	\omega_m &:=& \omega+b\label{eq:omega}\\
	\alpha_m &:=& \alpha\\
	\beta_m &:=& \beta.\label{eq:betam}
\end{IEEEeqnarray}

\section{The observer}\label{sec:Observer}
We show the state of the design system~\eqref{eq:alpha}--\eqref{eq:betam} can be estimated by the observer
\begin{IEEEeqnarray}{rCl}
	\dot{\hat\alpha} &=& \hat\alpha\times(\omega_m-\hat b)-k_\alpha(\hat\alpha-\alpha_m)\label{eq:alphahat}\\
	\dot{\hat\beta} &=& \hat\beta\times(\omega_m-\hat b)-k_\beta(\hat\beta-\beta_m)\\
	\dot{\hat b} &=& l_\alpha\hat\alpha\times\alpha_m+l_\beta\hat\beta\times\beta_m,\label{eq:bhat}
\end{IEEEeqnarray}
where $k_\alpha,k_\beta,l_\alpha,l_\beta>0$ are strictly positive constants. Notice this observer is very similar to the so-called bias observer of~\cite{BatisSO2012CEP}. The main difference is the use of the filtered terms $\hat\alpha\times(\omega_m-\hat b)$ and $\hat\beta\times(\omega_m-\hat b)$ instead of their unfiltered versions, which avoids injecting too much measurement noise when the gains are small; as a consequence the error system below is no longer linear in the measured quantities $\omega_m,\alpha_m,\beta_m$, and the technique of proof of~\cite{BatisSO2012CEP} is no longer applicable.

Defining the error vectors, $e_\alpha:=\hat\alpha-\alpha$, $e_\beta:=\hat\beta-\beta$ and $e_b:=\hat b-b$, the error system reads
\begin{IEEEeqnarray*}{rCl}
	\dot e_\alpha &=& e_\alpha\times\omega-(\alpha+e_\alpha)\times e_b-k_\alpha e_\alpha\label{eq:ealpha}\\
	\dot e_\beta &=& e_\beta\times\omega-(\beta+e_\beta)\times e_b-k_\beta e_\beta\label{eq:ebeta}\\
	\dot e_b &=& l_\alpha e_\alpha\times\alpha+l_\beta e_\beta\times\beta.\label{eq:eb}
\end{IEEEeqnarray*} 
It will be convenient to use the rotated variables $E_\alpha:=Re_\alpha$, $E_\beta:=Re_\beta$, $E_b:=Re_b$, and~$\Omega:=R\omega$. In terms of the rotated variables, the error system reads
\begin{IEEEeqnarray}{rCl}
	\dot E_\alpha &=& E_b\times(\alpha_i+E_\alpha)-k_\alpha E_\alpha\label{eq:Ealpha}\\
	\dot E_\beta &=& E_b\times(\beta_i+E_\beta)-k_\beta E_\beta\label{eq:Ebeta}\\
	\dot E_b &=& \Omega\times E_b+l_\alpha E_\alpha\times\alpha_i+l_\beta E_\beta\times\beta_i.\label{eq:Eb}
\end{IEEEeqnarray} 
Indeed, we find e.g. \eqref{eq:Ealpha} by writing
\begin{IEEEeqnarray*}{rCl}
	\dot E_\alpha &=& \dot Re_\alpha+R\dot e_\alpha\\
	&=& R(\omega\times e_\alpha)\\
	&& +\, R(e_\alpha\times\omega)-(R\alpha+Re_\alpha)\times Re_b-k_\alpha Re_\alpha\\
	&=& -(\alpha_i+E_\alpha)\times E_b-k_\alpha E_\alpha,
\end{IEEEeqnarray*} 
where we have used \eqref{eq:R}, \eqref{eq:ealpha} and $R(x\times y)=Rx\times Ry$ % whatever the vectors $x,y$ 
since $R$ is a rotation matrix.
It is of course equivalent to work with this rotated error system; notice also $\abs{\Omega}=\abs{\omega}$. 

\begin{thm}\label{th:conv2vec} Assume $k_\alpha,k_\beta,l_\alpha,l_\beta>0$ and $\omega$ bounded. % with an unknown bound. 
	Then the equilibrium point $(\bar E_\alpha,\bar E_\beta,\bar E_b):=(0,0,0)$ of the error system~\eqref{eq:Ealpha}--\eqref{eq:Eb} is uniformly globally asymptotically stable and locally exponentially stable.
\end{thm}

\begin{pf}	
Consider the candidate Lyapunov function 
\begin{IEEEeqnarray*}{rCl}
	V &:=& \sigma_1V_1+\sigma_2V_1^2+V_3\label{eq:LyapFct},
\end{IEEEeqnarray*}
where the coefficients $\sigma_1,\sigma_2>0$ are yet to be defined, and
\begin{IEEEeqnarray*}{rCl}
	V_1(E_\alpha,E_\beta,E_b) &:=& \frac{1}{2}\bigl(l_\alpha\abs{E_\alpha}^2 + l_\beta\abs{E_\beta}^2 + \abs{E_b}^2\bigr)\\
%	V_2(E_\alpha,E_\beta,E_b) &:=& V_1^2(E_\alpha,E_\beta,E_b)\\
	V_3(E_\alpha,E_\beta,E_b) &:=& \frac{1}{2}\abs{E_b - \frac{l_\alpha}{k_\alpha}\alpha_i\times E_\alpha - \frac{l_\beta}{k_\beta}\beta_i\times E_\beta}^2.
\end{IEEEeqnarray*} 	
Clearly, $V$ is positive definite and radially unbounded.

We now compute the time derivatives of its pieces along the trajectories of the error system~\eqref{eq:Ealpha}--\eqref{eq:Eb}.
\begin{IEEEeqnarray*}{rCl}
	\dot V_1&=& l_\alpha\inner{E_\alpha}{E_b\times\alpha_i+E_b\times E_\alpha - k_\alpha E_\alpha}\\
	&& +\> l_\beta\inner{E_\beta}{E_b\times\beta_i+E_b\times E_\beta - k_\beta E_\beta}\\
	&& +\> \inner{E_b}{\Omega\times E_b + l_\alpha E_\alpha\times\alpha_i + l_\beta E_\beta\times\beta_i}\\
	&=& -k_\alpha l_\alpha\abs{E_\alpha}^2 - k_\beta l_\beta\abs{E_\beta}^2,
\end{IEEEeqnarray*}
where we have used $\inner{x}{x\times y}=0$.
\begin{IEEEeqnarray*}{rCl}
	\frac{d}{dt}V_1^2 &=& -k_\alpha l_\alpha^2\abs{E_\alpha}^4 - k_\beta l_\beta^2\abs{E_\beta}^4
	- k_\alpha l_\alpha\abs{E_\alpha}^2\abs{E_b}^2\\
	&&-\> k_\beta l_\beta\abs{E_\beta}^2\abs{E_b}^2 - (k_\alpha+k_\beta)l_\alpha l_\beta\abs{E_\alpha}^2\abs{E_\beta}^2.
\end{IEEEeqnarray*}
Setting $u_{\alpha\beta}:=\frac{l_\alpha}{k_\alpha}\alpha_i\times E_\alpha + \frac{l_\beta}{k_\beta}\beta_i\times E_\beta$,
\begin{IEEEeqnarray*}{rCl}
	\IEEEeqnarraymulticol{3}{l}{ \inner{E_b}{\dot E_b-\dot u_{\alpha\beta}} }\\ ~~~
	&=& \inner{E_b}{\Omega\times E_b}+\inner{E_b}{l_\alpha E_\alpha\times\alpha_i+l_\beta E_\beta\times\beta_i}\\
	&& -\> \inner{E_b}{\frac{l_\alpha}{k_\alpha}\alpha_i\times(E_b\times\alpha_i) + \frac{l_\beta}{k_\beta}\beta_i\times(E_b\times\beta_i)}\\	
	&& -\> \inner{E_b}{\frac{l_\alpha}{k_\alpha}\alpha_i\times(E_b\times E_\alpha) + \frac{l_\beta}{k_\beta}\beta_i\times(E_b\times E_\beta)}\\
	&& +\> \inner{E_b}{\frac{l_\alpha}{k_\alpha}\alpha_i\times(k_\alpha E_\alpha) + \frac{l_\beta}{k_\beta}\beta_i\times(k_\beta E_\beta)}\\
	&=& \inner{E_b}{\Mab E_b} \\
	&& -\> \inner{E_b}{\frac{l_\alpha}{k_\alpha}\alpha_i\times(E_b\times E_\alpha)
		+\frac{l_\beta}{k_\beta}\beta_i\times(E_b\times E_\beta)},
\end{IEEEeqnarray*}
where we have used $\inner{x}{x\times y}=0$ and $\inner{x}{y\times z}=\inner{z}{x\times y}$.
\begin{IEEEeqnarray*}{rCl}
	\IEEEeqnarraymulticol{3}{l}{ \inner{u_{\alpha\beta}}{\dot E_b-\dot u_{\alpha\beta}} }\\ ~~~
	&=& \inner{u_{\alpha\beta}}{\Omega\times E_b} + \inner{u_{\alpha\beta}}{\Mab E_b} \\
	&& -\> \inner{u_{\alpha\beta}}{\frac{l_\alpha}{k_\alpha}\alpha_i\times(E_b\times E_\alpha)
		+\frac{l_\beta}{k_\beta}\beta_i\times(E_b\times E_\beta)}\\
	&=& \inner{E_b}{u_{\alpha\beta}\times\Omega} + \inner{E_b}{\Mab u_{\alpha\beta}} \\
	&& -\> \inner{E_b}{\frac{l_\alpha}{k_\alpha}\alpha_i\times(u_{\alpha\beta}\times E_\alpha)
	+\frac{l_\beta}{k_\beta}\beta_i\times(u_{\alpha\beta}\times E_\beta)},
\end{IEEEeqnarray*}
where we have used $\inner{x}{y\times z}=\inner{z}{x\times y}$ several times.

We next bound the last two expressions. Since $\alpha_i,\beta_i$ are independent and $k_\alpha,l_\alpha,k_\beta,l_\beta>0$, the matrix %$-\Mab$
\begin{IEEEeqnarray*}{C}
	-\Mab
\end{IEEEeqnarray*}
is positive definite, see~\cite[Lemma~2]{TayebRB2013TAC}. As a consequence, there exists $\mu>0$ such that
\begin{IEEEeqnarray*}{rCl}
	\inner{E_b}{\Bigl(\frac{l_\alpha}{k_\alpha}a_{i\times}^2 +\frac{l_\beta}{k_\beta}b_{i\times}^2\Bigr)E_b}
	&\le& -\mu\abs{E_b}^2.
\end{IEEEeqnarray*}

Using then Young's inequality $\abs{\inner{x}{y}}\leq \frac{\epsilon\abs{x}^2}{2}+\frac{\abs{y}^2}{2\epsilon}$ and
$\bigabs{\sum_{i=1}^nx_i}^2\le n\sum_{i=1}^n\abs{x_i}^2$ yields
\begin{IEEEeqnarray*}{rCl}
	\IEEEeqnarraymulticol{3}{l}{ \inner{E_b}{\dot E_b-\dot u_{\alpha\beta}} }\\ ~
	&\le& -\mu\abs{E_b}^2 + \frac{\varepsilon}{2}\abs{E_b}^2\\
	&& +\>\frac{1}{\varepsilon}\Bigabs{ \frac{l_\alpha}{k_\alpha}\alpha_i\times(E_b\times E_\alpha) }^2
	+\frac{1}{\varepsilon}\Bigabs{ \frac{l_\beta}{k_\beta}\beta_i\times(E_b\times E_\beta) }^2\\
	&\le& \Bigl(\frac{\varepsilon}{2}-\mu\Bigr)\abs{E_b}^2\\
	&& +\> \frac{l_\alpha^2\abs{\alpha_i}^2}{\varepsilon k_\alpha^2}\abs{E_\alpha}^2\abs{E_b}^2
	+\frac{l_\beta^2\abs{\beta_i}^2}{\varepsilon k_\beta^2}\abs{E_\beta}^2\abs{E_b}^2.
\end{IEEEeqnarray*}
Similarly,
\begin{IEEEeqnarray*}{rCl}
	\IEEEeqnarraymulticol{3}{l}{ \bigabs{\inner{u_{\alpha\beta}}{\dot E_b-\dot u_{\alpha\beta}}} }\\ ~
	&\le& \frac{\varepsilon}{2}\abs{E_b}^2 + \frac{2}{\varepsilon}\abs{ u_{\alpha\beta}\times\Omega }^2
	+ \frac{2}{\varepsilon}\Bigabs{ \Mab u_{\alpha\beta} }^2\\
	&& +\> \frac{2}{\varepsilon}\Bigabs{ \frac{l_\alpha}{k_\alpha}\alpha_i\times(u_{\alpha\beta}\times E_\alpha) }^2
	+\frac{2}{\varepsilon}\Bigabs{ \frac{l_\beta}{k_\beta}\beta_i\times(u_{\alpha\beta}\times E_\beta) }^2\\
	&\le& \frac{\varepsilon}{2}\abs{E_b}^2 + \frac{2c_\omega^2}{\varepsilon}\abs{u_{\alpha\beta}}^2
	+ \frac{4}{\varepsilon}\Bigl(\frac{l_\alpha^2}{k_\alpha^2}\abs{\alpha_i}^4 + \frac{l_\beta^2}{k_\beta^2}\abs{\beta_i}^4\Bigr)\abs{ u_{\alpha\beta} }^2\\
	&& +\> \frac{2l_\alpha^2\abs{\alpha_i}^2}{\varepsilon k_\alpha^2}\abs{ E_\alpha }^2\abs{ u_{\alpha\beta} }^2
	+ \frac{2l_\beta^2\abs{\beta_i}^2}{\varepsilon k_\beta^2}\abs{ E_\beta }^2\abs{ u_{\alpha\beta} }^2,
\end{IEEEeqnarray*}
where $\abs{\Omega}=\abs{\omega}\le c_\omega$ and
\begin{IEEEeqnarray*}{rCl}
	\abs{u_{\alpha\beta}}^2 &\le& \frac{2l_\alpha^2\abs{\alpha_i}^2}{k_\alpha^2}\abs{ E_\alpha }^2
	+ \frac{2l_\beta^2\abs{\beta_i}^2}{k_\beta^2}\abs{ E_\beta }^2.
\end{IEEEeqnarray*}

Collecting all the pieces, we eventually find
\begin{IEEEeqnarray*}{rCl}
	\dot V &\le& -\mu'\abs{E_b}^2 - \sigma_{1\alpha}\abs{ E_\alpha }^2 - \sigma_{1\beta}\abs{ E_\beta }^2\\
	&&-\> \sigma_{2\alpha}\abs{ E_\alpha }^4 - \sigma_{2\beta}\abs{ E_\beta }^4 - \sigma_{2\alpha\beta}\abs{ E_\alpha }^2\abs{ E_\beta }^2\\
	&&-\> \sigma'_{2\alpha}\abs{E_\alpha}^2\abs{E_b}^2 - \sigma'_{2\beta}\abs{E_\beta}^2\abs{E_b}^2,
\end{IEEEeqnarray*}
where
\begin{IEEEeqnarray*}{rCl}
	\mu' &:=& \mu-\varepsilon\\
	\sigma_{1\alpha} &:=& \sigma_1k_\alpha l_\alpha - \frac{4}{\varepsilon}\Bigl(c_\omega^2+\frac{2l_\alpha^2}{k_\alpha^2}\abs{\alpha_i}^4 
	+ \frac{2l_\beta^2}{k_\beta^2}\abs{\beta_i}^4\Bigr)\frac{l_\alpha^2\abs{\alpha_i}^2}{k_\alpha^2}\\
	\sigma_{1\beta} &:=& \sigma_1k_\beta l_\beta - \frac{4}{\varepsilon}\Bigl(c_\omega^2+\frac{2l_\alpha^2}{k_\alpha^2}\abs{\alpha_i}^4 
	+ \frac{2l_\beta^2}{k_\beta^2}\abs{\beta_i}^4\Bigr)\frac{l_\beta^2\abs{\beta_i}^2}{k_\beta^2}\\
	\sigma_{2\alpha} &:=& \sigma_2k_\alpha l_\alpha^2 - \frac{4l_\alpha^4\abs{\alpha_i}^4}{\varepsilon k_\alpha^4}\\
	\sigma_{2\beta} &:=& \sigma_2k_\beta l_\beta^2 - \frac{4l_\beta^4\abs{\beta_i}^4}{\varepsilon k_\beta^4}\\
	\sigma_{2\alpha\beta} &:=& \sigma_2(k_\alpha+k_\beta)l_\alpha l_\beta 
	- \frac{8l_\alpha^2l_\beta^2\abs{\alpha_i}^2\abs{\beta_i}^2}{\varepsilon k_\alpha^2k_\beta^2}\\
	\sigma'_{2\alpha} &:=& \sigma_2k_\alpha l_\alpha - \frac{l_\alpha^2\abs{\alpha_i}^2}{\varepsilon k_\alpha^2}\\
	\sigma'_{2\beta} &:=& \sigma_2k_\beta l_\beta - \frac{l_\beta^2\abs{\beta_i}^2}{\varepsilon k_\beta^2}.
\end{IEEEeqnarray*}
All the above coefficients are strictly positive if $\varepsilon$ is chosen small enough and $\sigma_1,\sigma_2$ large enough. This means $V$ is a strict Lyapunov function, which proves uniform global asymptotic stability. Notice the bound $c_\omega$ need not be known, since $\sigma_1$ can always been chosen large enough to achieve $\sigma_{1\alpha},\sigma_{1\beta}>0$

To establish local exponential stability, we consider the linearized error system around the equilibrium point $(\bar E_\alpha,\bar E_\beta,\bar E_b):=(0,0,0)$; it reads
\begin{IEEEeqnarray*}{rCl}
	\delta\dot E_\alpha &=& \delta E_b\times\alpha_i-k_\alpha\delta E_\alpha\\
	\delta\dot E_\beta &=& \delta E_b\times\beta_i-k_\beta E_\beta\\
	\delta\dot E_b &=& \Omega\times\delta E_b+l_\alpha\delta E_\alpha\times\alpha_i+l_\beta\delta E_\beta\times\beta_i.
\end{IEEEeqnarray*} 
This is a linear time-varying system (due to the dependence on $\Omega$); its origin is uniformly locally asymptotically stable, since the origin of the nonlinear system is uniformly globally asymptotically stable. Exponential stability follows from the fact that for linear (time-varying) systems, uniform asymptotic stability and exponential stability are equivalent, see~\cite[Theorem 4.11]{khalil2002book}.\qed
\end{pf}

\begin{rem}
	$V_1$ alone is a (non strict) Lyapunov function for the error system~\eqref{eq:Ealpha}--\eqref{eq:Eb}; using repeatedly Barbalat's lemma, (non uniform) global asymptotic stability can be very easily established, see~\cite{MartiS2016arXiv}. On the other hand, the more complicated Lyapunov function $V$ used in the proof is strict, hence yields uniform stability.
\end{rem} 

\begin{rem}
	If in the observer~\eqref{eq:alphahat}-\eqref{eq:bhat} the filtered terms $\hat\alpha\times(\omega_m-\hat b)$ and $\hat\beta\times(\omega_m-\hat b)$ are replaced by their filtered versions $\hat\alpha\times\omega_m-\alpha_m\times\hat b$ and $\hat\beta\times\omega_m-\beta_m\times\hat b$, we recover the first part of the linear cascaded observer of~\cite{BatisSO2012CEP}. In this case, the strict Lyapunov function $V$ used in the proof can also be used to establish uniform global asymptotic stability, which provides a simpler and more constructive proof than the ``abstract'' proof of~\cite{BatisSO2012CEP} based on uniform observability arguments.
\end{rem} 

\begin{rem}
	The result obviously also holds if the scalar gains $k_\alpha,k_\beta,l_\alpha,l_\beta$ are replaced by $3\times3$ symmetric definite positive matrices; these added degrees of freedom for the tuning might be useful in practice if the components of the measurement vectors are produced by sensors with very different characteristics.
	
	It is also clear that more than two vectors can be used with an obvious generalization of the proposed structure.
\end{rem}

\begin{rem}
	The observer does not use the knowledge of the constant vectors $\alpha_i$ and $\beta_i$. This may be an interesting feature in some applications when those vectors for example are not precisely known and/or (slowly) vary.
\end{rem}

We then have the obvious but important following corollary, which gives an estimate of the true orientation matrix~$R$ by using the knowledge of $\alpha_i$ and $\beta_i$. Notice it is considerably simpler than the approach proposed in~\cite{BatisSO2012CEP}, where the estimated orientation matrix is obtained through an additional observer of dimension~$9$.
\begin{cor} Under the assumptions of theorem~\ref{th:conv2vec}, the matrix $\tilde R$ defined by
	\begin{IEEEeqnarray*}{rCl}
		\tilde R^T &:=& \begin{pmatrix}\frac{\hat\alpha}{\abs{\alpha_i}}& \frac{\hat\alpha\times\hat\beta}{\abs{\alpha_i\times\beta_i}} & \frac{\hat\alpha\times(\hat\alpha\times\hat\beta)}{\abs{\alpha_i\times(\alpha_i\times\beta_i)}}\end{pmatrix}\cdot R_i^T\\
		R_i &:=& \begin{pmatrix}\frac{\alpha_i}{\abs{\alpha_i}}& \frac{\alpha_i\times\beta_i}{\abs{\alpha_i\times\beta_i}} & \frac{\alpha_i\times(\alpha_i\times\beta_i)}{\abs{\alpha_i\times(\alpha_i\times\beta_i)}}\end{pmatrix}
	\end{IEEEeqnarray*}	
	uniformly globally asymptotically converges to~$R$.
\end{cor}
\begin{pf}
	By Theorem~\ref{th:conv2vec}, $\hat\alpha\to\alpha$ and $\hat\beta\to\beta$. Hence,
	\begin{IEEEeqnarray*}{rCl}
		\tilde R^T &\to& \begin{pmatrix}\frac{\alpha}{\abs{\alpha_i}}& \frac{\alpha\times\beta}{\abs{\alpha_i\times\beta_i}} & \frac{\alpha\times(\alpha\times\beta)}{\abs{\alpha_i\times(\alpha_i\times\beta_i)}}\end{pmatrix}\cdot R_i^T\\
		&=& \begin{pmatrix}\frac{R^T\alpha_i}{\abs{\alpha_i}}& \frac{R^T\alpha_i\times R^T\beta_i}{\abs{\alpha_i\times\beta_i}} & \frac{R^T\alpha_i\times(R^T\alpha_i\times R^T\beta_i)}{\abs{\alpha_i\times(\alpha_i\times\beta_i)}}\end{pmatrix}\cdot R_i^T\\
		&=& R^TR_iR_i^T\\
		&=& R^T,
	\end{IEEEeqnarray*}	
	where we have used $R^T(u\times v)=R^Tu\times R^Tv$ since $R$ is a rotation matrix.\qed
\end{pf}

Of course, $\tilde R^T$ has no reason to be a rotation matrix (it is only asymptotically so); it is nevertheless the product of a matrix with orthogonal (possibly zero) columns by a rotation matrix. If a bona fide rotation matrix is required at all times, a natural idea is to project $\tilde R$ on the ``closest'' rotation matrix~$\hat R$, thanks to a polar decomposition. Because of the particular form of $\tilde R$, the expression of $\hat R$ can readily be found without using the standard but computationally heavy projection algorithm based on singular value decomposition. For details about the polar decomposition and related matters, see e.g.~\cite[Chapter~$8$]{Higham2008book}.
\begin{prop}
	Consider the polar decomposition of $\tilde R^T$
	\begin{IEEEeqnarray*}{rCl}
		\tilde R^T &=& \hat R^T(\tilde R\tilde R^T)^\frac{1}{2}.
	\end{IEEEeqnarray*} 
	$\hat R$, which is by construction the best approximation of~$\tilde R$ among all orthogonal matrices, is a rotation matrix that uniformly globally asymptotically converges to~$R$. When $\hat\alpha$ and $\hat\beta$ are not collinear, $\hat R$ is uniquely defined by
	\begin{IEEEeqnarray*}{rCl}
		\hat R^T &:=& \begin{pmatrix}\frac{\hat\alpha}{\abs{\hat\alpha}}& \frac{\hat\alpha\times\hat\beta}{\abs{\hat\alpha\times\hat\beta}} & \frac{\hat\alpha\times(\hat\alpha\times\hat\beta)}{\abs{\hat\alpha\times(\hat\alpha\times\hat\beta)}}\end{pmatrix}\cdot R_i^T.
	\end{IEEEeqnarray*}	
\end{prop}
\begin{pf}
	Since $\tilde R$ is the product of a matrix with orthogonal columns by a rotation matrix,
	\begin{IEEEeqnarray*}{rCl}
		(\tilde R\tilde R^T)^\frac{1}{2}
		&=& \begin{pmatrix}R_i
		\begin{pmatrix}\frac{\abs{\hat\alpha}^2}{\abs{\alpha_i}^2}& 0&0\\
			0& \frac{\abs{\hat\alpha\times\hat\beta}^2}{\abs{\alpha_i\times\beta_i}^2} & 0\\
			0&0& \frac{\abs{\hat\alpha\times(\hat\alpha\times\hat\beta)}^2}{\abs{\alpha_i\times(\alpha_i\times\beta_i)}^2}\end{pmatrix}
		R_i^T\end{pmatrix}^\frac{1}{2}\\
		&=& R_i
	\begin{pmatrix}\frac{\abs{\hat\alpha}}{\abs{\alpha_i}}& 0&0\\
		0& \frac{\abs{\hat\alpha\times\hat\beta}}{\abs{\alpha_i\times\beta_i}} & 0\\
		0&0& \frac{\abs{\hat\alpha\times(\hat\alpha\times\hat\beta)}}{\abs{\alpha_i\times(\alpha_i\times\beta_i)}}\end{pmatrix}
		R_i^T.		
	\end{IEEEeqnarray*} 
	When $\hat\alpha$ and $\hat\beta$ are not collinear, the expression for $\hat R$ follows at once from $\hat R^T=\tilde R^T(\tilde R\tilde R^T)^{-\frac{1}{2}}$. When $\hat\alpha=0$, one may choose $\hat R^T:=R_i^T$; when $\hat\alpha\neq0$ but $\hat\alpha\times\hat\beta=0$, one may choose $\hat R^T:=(\frac{\hat\alpha}{\abs{\hat\alpha}},\hat E_2,\hat E_3)\cdot R_i^T$, where $\frac{\hat\alpha}{\abs{\hat\alpha}}$, $\hat E_2$ and $\hat E_3$ form a direct orthonormal frame.\qed
\end{pf}

\begin{figure}[ht]\hspace{-5mm}
	\includegraphics[width=1.1\columnwidth]{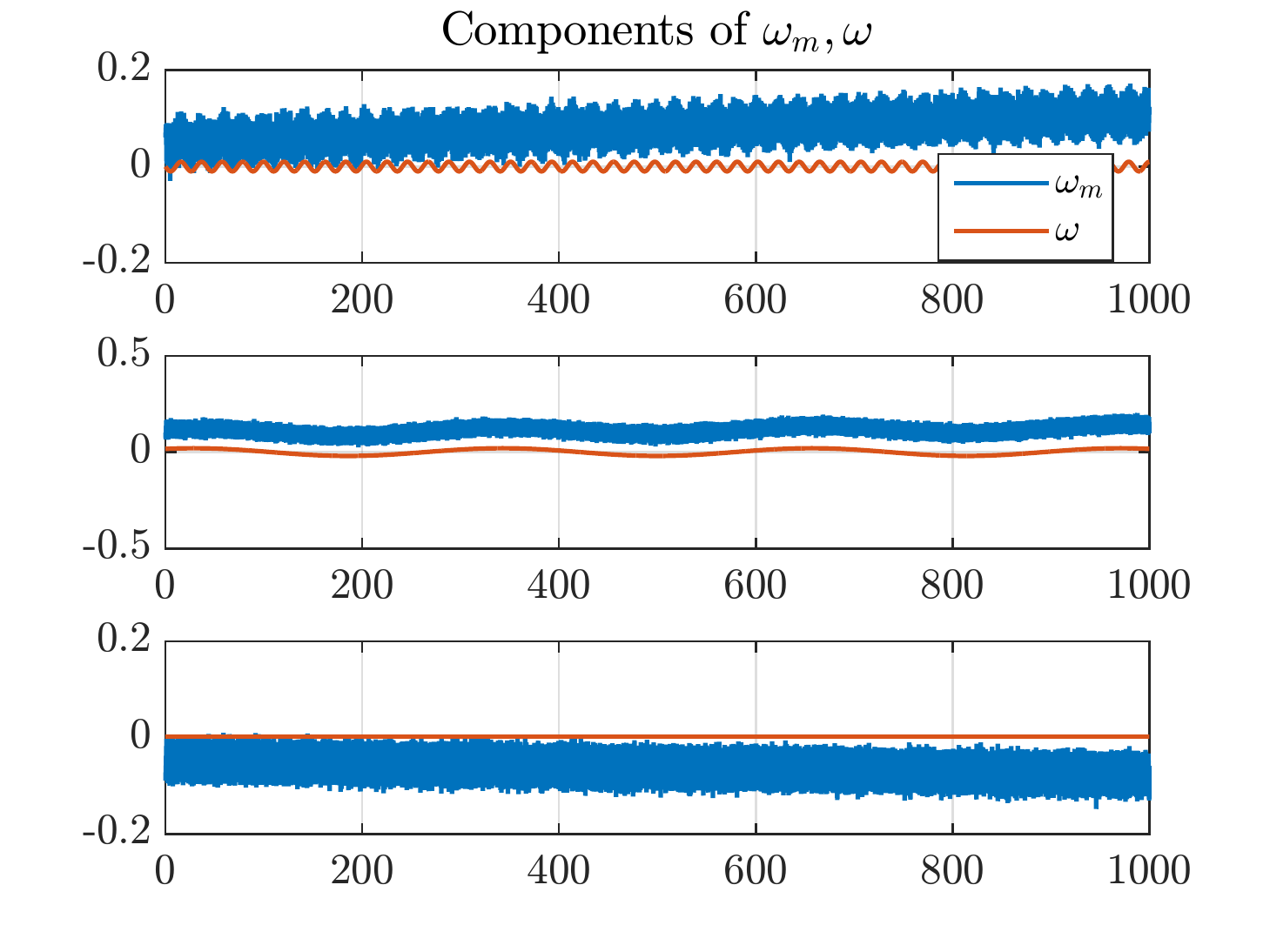}
	\caption{Components of true $\omega$ (red) and measured $\omega_m$ (blue).}
	\label{fig:omega}
\end{figure}

\begin{figure}[ht]\hspace{-5mm}
	\includegraphics[width=1.1\columnwidth]{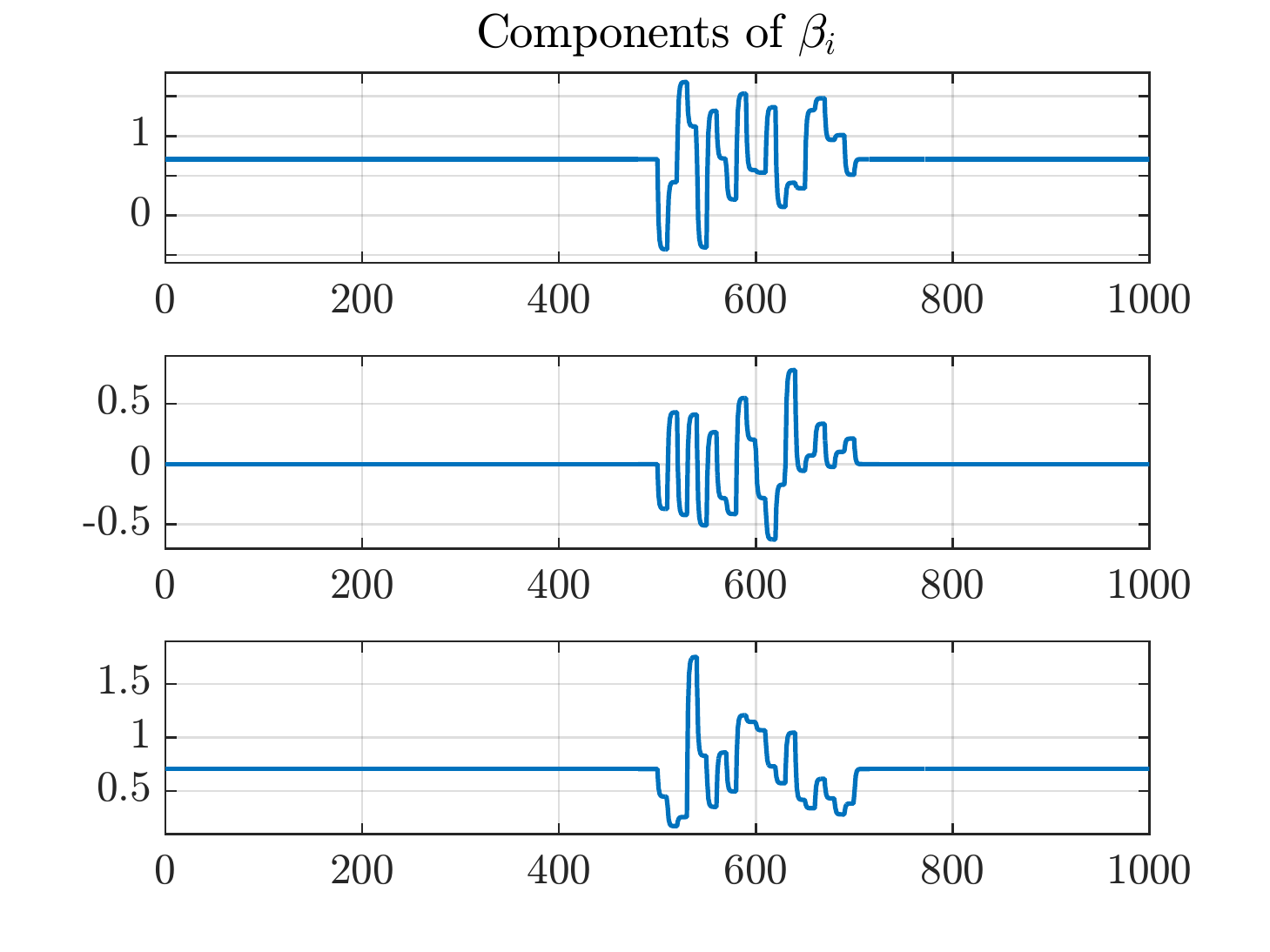}
	\caption{Components of ``constant'' vector~$\beta_i$.}
	\label{fig:betai}
\end{figure}

\begin{figure}[ht]\hspace{-5mm}
	\includegraphics[width=1.1\columnwidth]{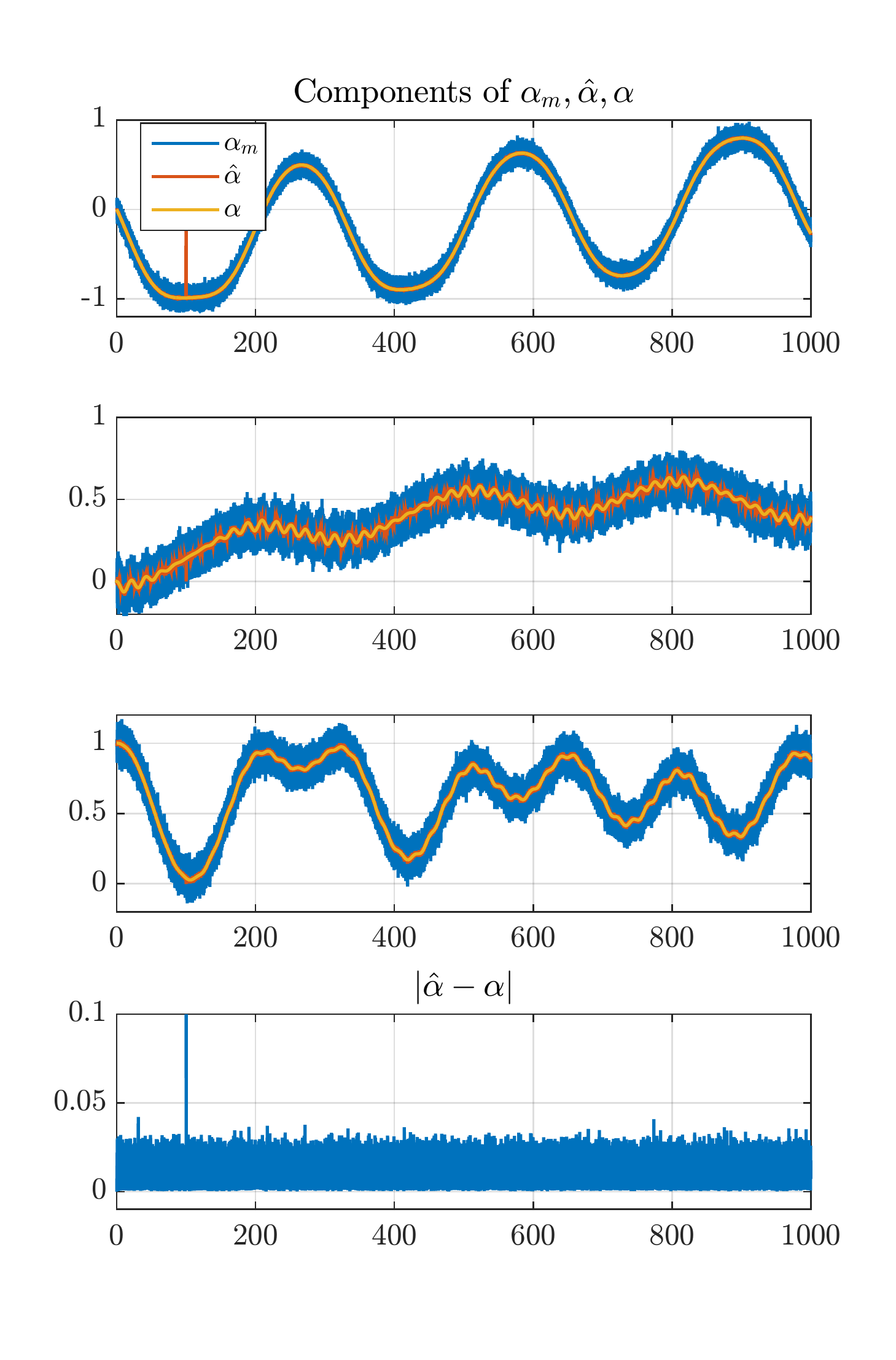}
	\caption{Components of true $\alpha$ (red), measured $\alpha_m$ (blue) and estimated~$\hat\alpha$ (orange).}
	\label{fig:alpha}
\end{figure}

\begin{figure}[ht]\hspace{-5mm}
	\includegraphics[width=1.1\columnwidth]{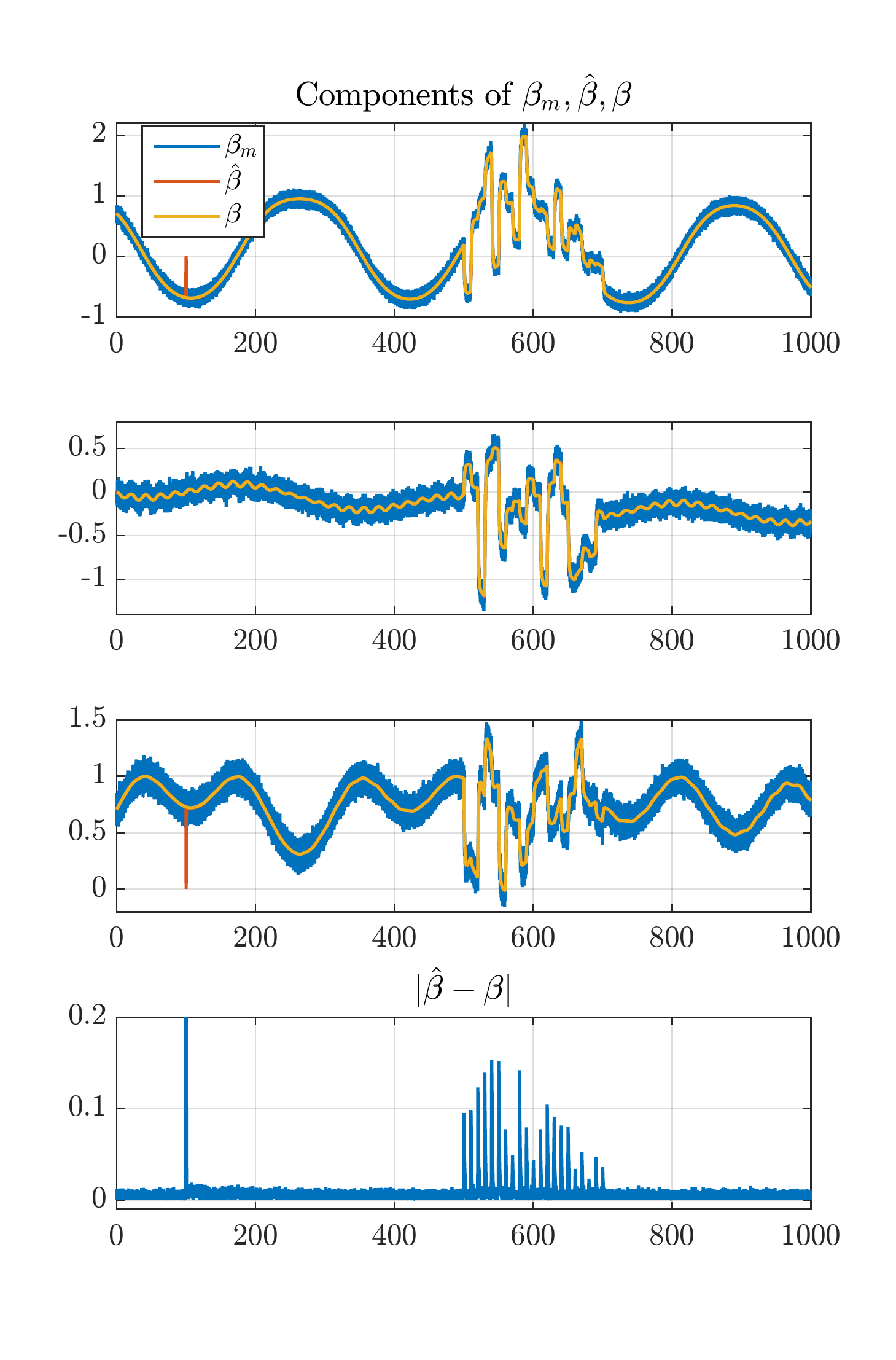}
	\caption{Components of true $\beta$ (red), measured $\beta_m$ (blue) and estimated~$\hat\beta$ (orange).}
	\label{fig:beta}
\end{figure}

\section{Simulations}\label{sec:Simulations}
The good behavior of the observer is now illustrated in simulation. The system starts in the initial state $R(0):=\mathrm{I}$, i.e. $\bigl(\phi(0),\theta(0),\psi(0)\bigr):=(0,0,0)$, and then undergoes the angular velocity~$\omega(t)$ displayed in Fig.~\ref{fig:omega}. The constant vectors $\alpha_i$ and $\beta_i$ are  respectively  set to the nominal values $(0,0,1)^T$ and $(\frac{1}{\sqrt2},0,\frac{1}{\sqrt2})^T$, which mimics the gravity and magnetic vectors. Moreover, $\beta_i$ is subjected to violent disturbances for $t\in[500,700]$, see Fig.~\ref{fig:betai}; of course only the nominal values of $\alpha_i$ and $\beta_i$ are known to the observer.

\begin{figure}[ht]\hspace{-5mm}
	\includegraphics[width=1.1\columnwidth]{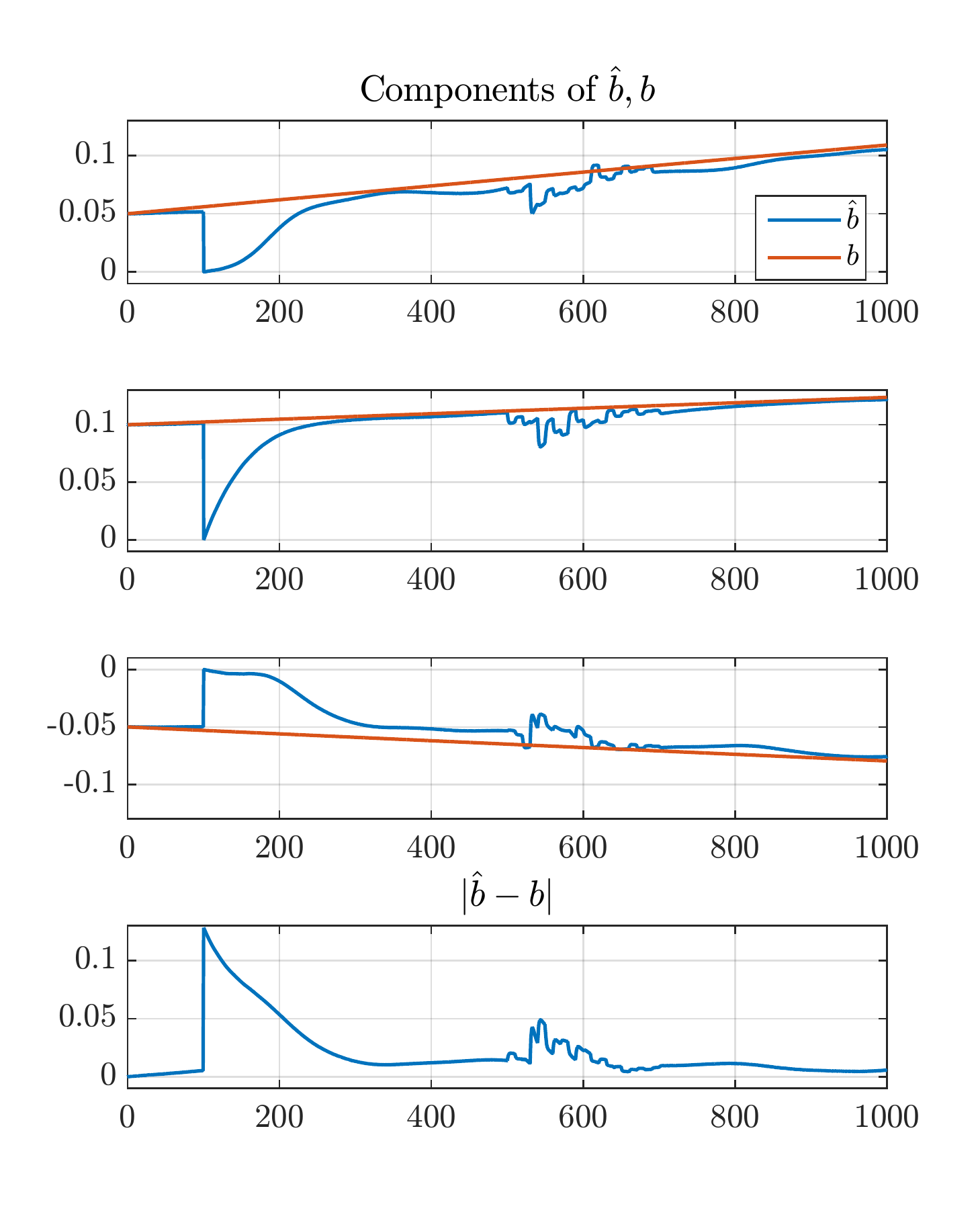}
	\caption{Components of true $b$ (red) and estimated $\hat b$ (blue).}
	\label{fig:biases}
\end{figure}

The observer is fed with the measured signals $\omega_m,\alpha_m$ and $\beta_m$, see Fig.~\ref{fig:omega}-\ref{fig:alpha}-\ref{fig:beta}; the measured velocity is affected by the (unknown) slowly drifting bias~$b$, see Fig.~\ref{fig:omega}-\ref{fig:biases}. All the measurement signals are corrupted by band-limited independent gaussian white noises with sample time \num{e-3} and rather large noise powers (\num{2e-6} for the components of~$\alpha_m,\beta_m$, and \num{2e-7} for those of~$\omega_m$). The tuning gains are set to $(k_\alpha,k_\beta,l_\alpha,l_\beta)=(10,10,0.15,0.15)$.

The observer is initialized with no error, but suddenly reinitialized to zero at~$t=100$. The convergence of the estimated state $\hat\alpha,\hat\beta$ and $\hat b$ is as anticipated excellent after the reinitialization, very fast for $\hat\alpha,\hat\beta$ and slower for $\hat b$, in accordance with the choice of gains, see Fig.~\ref{fig:alpha}-\ref{fig:beta}-\ref{fig:biases}. The convergence is also very good for $t\in[500,700]$ when $\beta_i$ is violently disturbed: the error $\hat\beta-\beta$ exhibits only small spikes, $\hat b$ is only slightly affected, and $\hat\alpha$ not affected at all at this scale; this illustrates the (desirable) independence between $\hat\alpha$ and $\hat\beta$, which are only slightly coupled through~$\hat b$. Notice also that the disturbances on $\beta_i$ are interpreted by the observer as a variation of~$b$.

\begin{figure}[ht]\hspace{-5mm}
	\includegraphics[width=1.1\columnwidth]{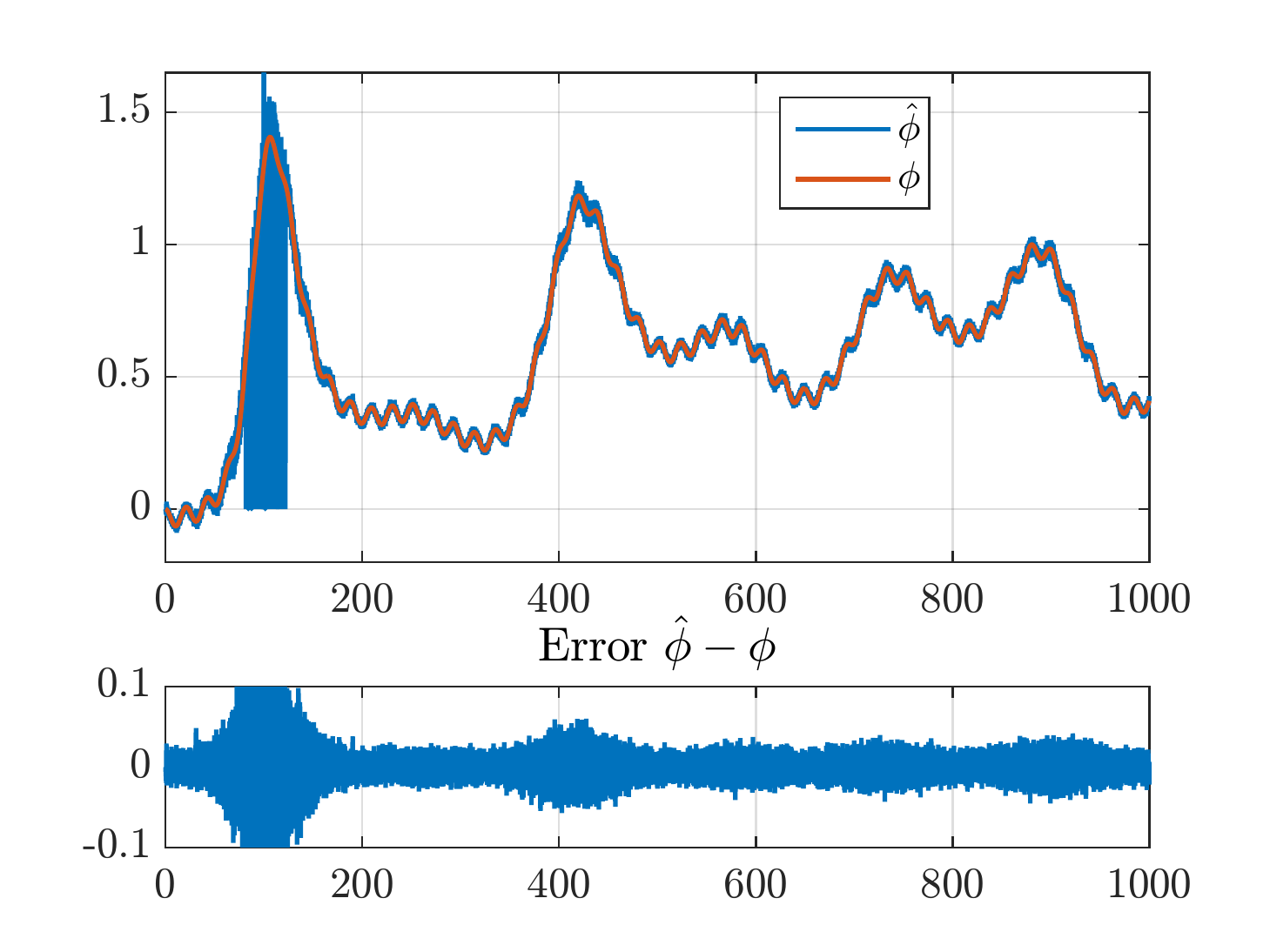}
	\caption{True $\phi$ (red) and estimated $\hat\phi$ (blue).}
	\label{fig:phi}
\end{figure}

\begin{figure}[ht]\hspace{-5mm}
	\includegraphics[width=1.1\columnwidth]{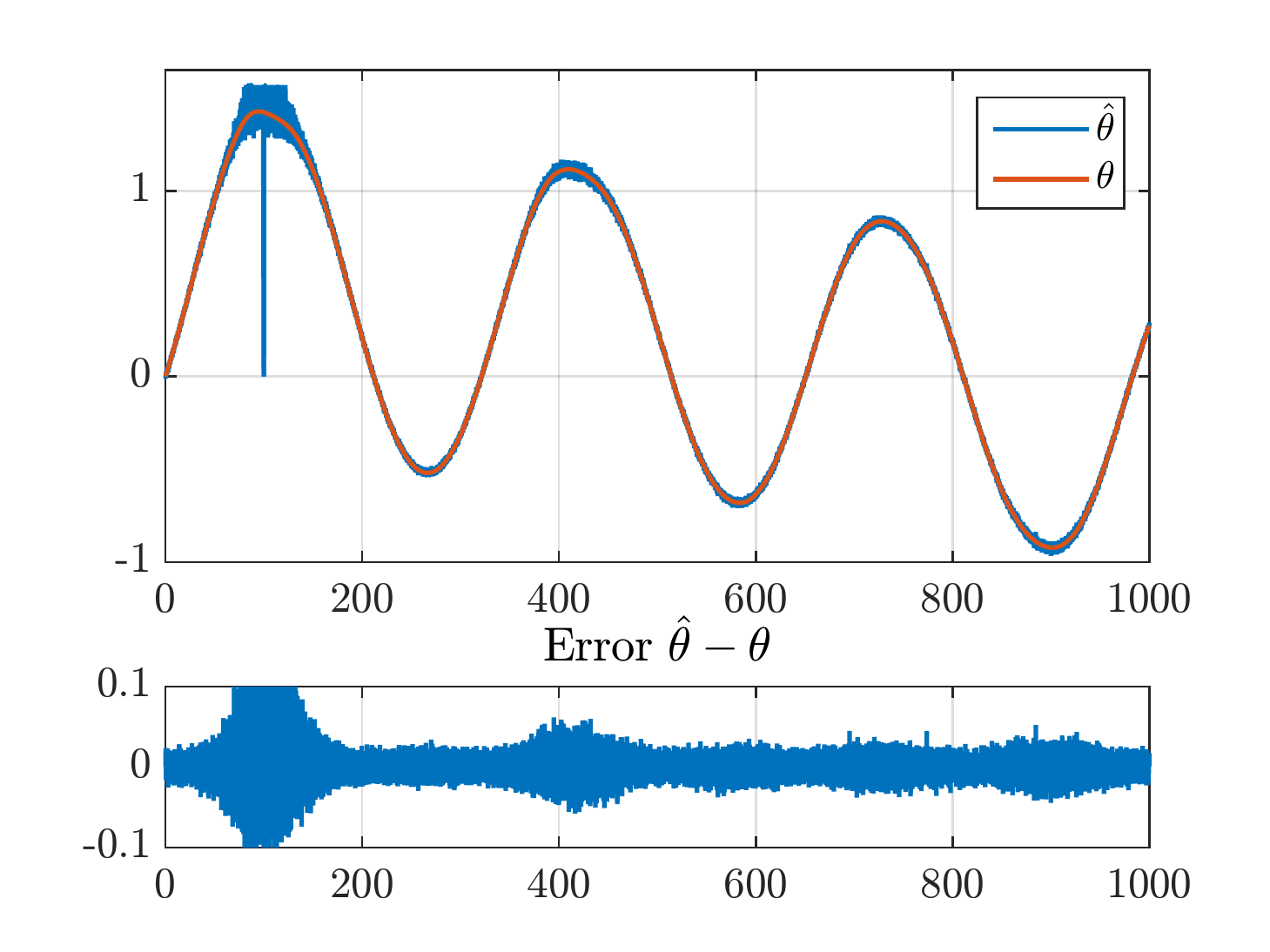}
	\caption{True $\theta$ (red) and estimated $\hat\theta$ (blue).}
	\label{fig:theta}
\end{figure}

\begin{figure}[ht]\hspace{-5mm}
	\includegraphics[width=1.1\columnwidth]{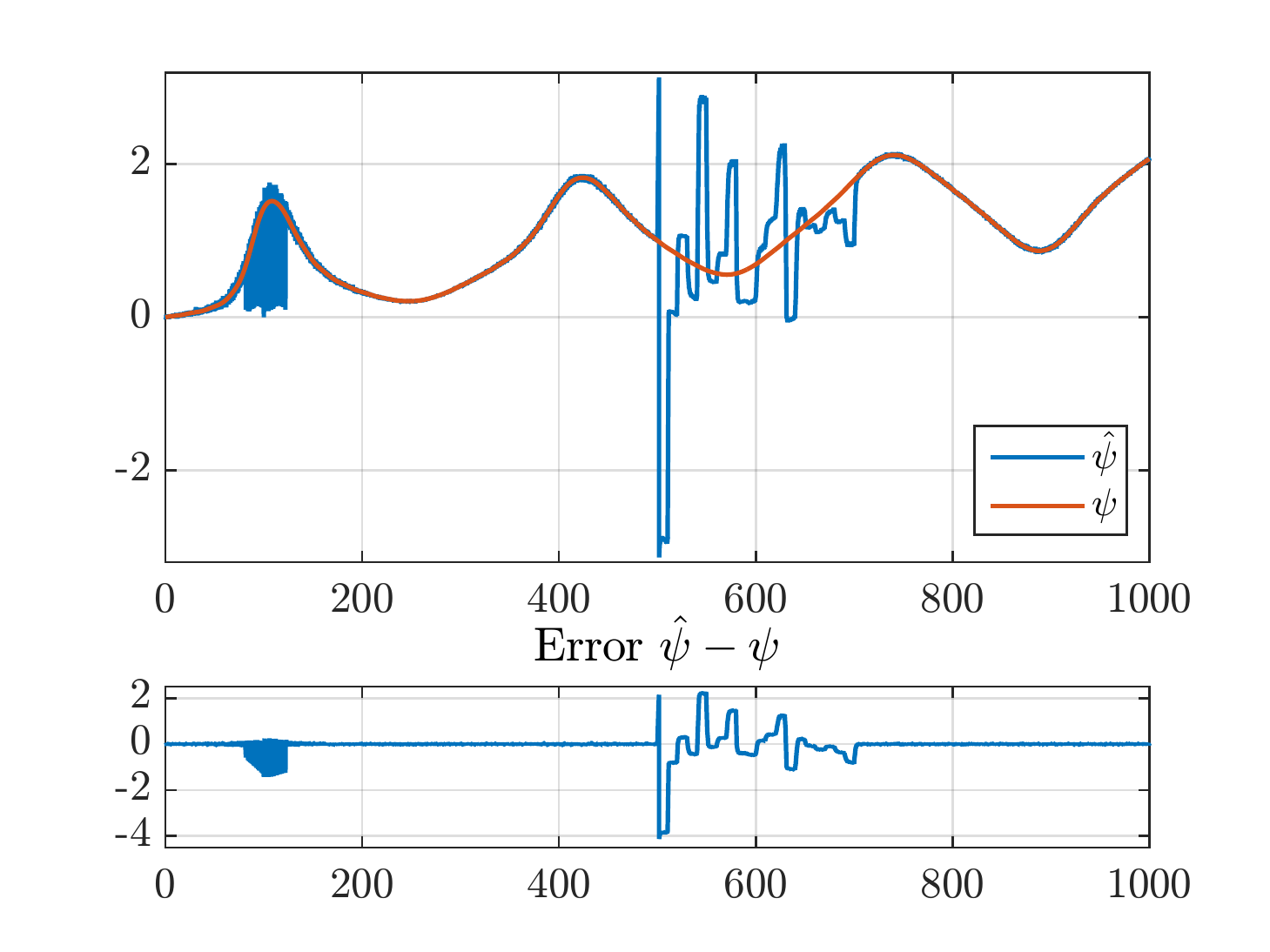}
	\caption{True $\psi$ (red) and estimated $\hat\phi$ (blue).}
	\label{fig:psi}
\end{figure}

We insist that the observer does not need the knowledge of the ``constant'' vectors $\alpha_i$ and $\beta_i$, which is why it can converge to the true state even when $\alpha_i,\beta_i$ vary not too fast (with respect to the gains $k_\alpha,k_\beta$). These vectors are needed only when the rotation matrix $R$ (or Euler angles, or quaternions) must be reconstructed. Fig.~\ref{fig:phi}-\ref{fig:theta}-\ref{fig:psi} show the reconstruction of the Euler angles $\phi,\theta,\psi$ (in radians) from $\hat\alpha,\hat\beta$ using the nominal values of $\alpha_i,\beta_i$. Notice the pitch angle $\theta$ and roll angle~$\phi$ (which are paramount for the control of for instance UAVs) are perfectly estimated even during the period where $\beta_i$ is violently disturbed (the fast oscillations of $\hat\theta$ just after $t=100$ are due to the angle determinacy which switches between $0$ and $\pi$, and have nothing to do with the observer); inevitably the estimated yaw angle~$\hat\phi$ is severely affected, since a completely wrong value of $\beta_i$ is used. This scenario illustrates an interesting practical feature of the observer, not enjoyed by e.g. the observer in~\cite{MahonHP2008TAC}: even if the magnetic sensor (i.e., $\beta_m$) is easily disturbed, it is nevertheless worth using to help the accelerometer (i.e., $\alpha_m$) in the estimation of the gyro biases, hence in the reconstruction of the pitch and roll angles, whatever the excitation (or absence of) provide by the angular velocity.

\section{Conclusion}\label{sec:conclusions}
We have presented a simple nonlinear ``geometry-free'' observer for attitude and gyro bias estimation, with guaranteed uniform global convergence and local exponential convergence. Simulations demonstrate that it performs very well, even with noisy measurements and not so constant inertial vectors $\alpha_i$ and~$\beta_i$ and bias~$b$. It can be seen as an interesting alternative to the $\SO$-based observer of~\cite{MahonHP2008TAC} or of the more complicated ``geometry-free'' observer of~\cite{BatisSO2012CEP}.

\ack
We thank our colleague Laurent Praly for enlightening discussions on the construction of strict Lyapunov functions. 
%________________________________________________________
\bibliography{ref}

\begin{thebibliography}{16}
\providecommand{\natexlab}[1]{#1}
\providecommand{\url}[1]{\texttt{#1}}
\providecommand{\urlprefix}{URL }
\expandafter\ifx\csname urlstyle\endcsname\relax
  \providecommand{\doi}[1]{doi:\discretionary{}{}{}#1}\else
  \providecommand{\doi}{doi:\discretionary{}{}{}\begingroup
  \urlstyle{rm}\Url}\fi

\bibitem[{Batista et~al.(2012{\natexlab{a}})Batista, Silvestre, and
  Oliveira}]{BatisSO2014AUT}
Batista, P., Silvestre, C., and Oliveira, P. (2012{\natexlab{a}}).
\newblock A {GES} attitude observer with single vector observations.
\newblock \emph{Automatica}, 48(2), 388--395.

\bibitem[{Batista et~al.(2012{\natexlab{b}})Batista, Silvestre, and
  Oliveira}]{BatisSO2012CEP}
Batista, P., Silvestre, C., and Oliveira, P. (2012{\natexlab{b}}).
\newblock Globally exponentially stable cascade observers for attitude
  estimation.
\newblock \emph{Control Engineering Practice}, 20(2), 148 -- 155.

\bibitem[{Bonnabel et~al.(2008)Bonnabel, Martin, and Rouchon}]{BoMaRo_tac09}
Bonnabel, S., Martin, P., and Rouchon, P. (2008).
\newblock Symmetry-preserving observers.
\newblock \emph{IEEE Transactions on Automatic Control}, 53(11), 2514--2526.

\bibitem[{Crassidis et~al.(2007)Crassidis, Markley, and
  Cheng}]{CrassMC2007JGCD}
Crassidis, J.L., Markley, F.L., and Cheng, Y. (2007).
\newblock Survey of nonlinear attitude estimation methods.
\newblock \emph{Journal of Guidance, Control, and Dynamics}, 30(1), 12--28.

\bibitem[{Grip et~al.(2012)Grip, Fossen, Johansen, and Saberi}]{GripFJS2012TAC}
Grip, H.F., Fossen, T.I., Johansen, T.A., and Saberi, A. (2012).
\newblock Attitude estimation using biased gyro and vector measurements with
  time-varying reference vectors.
\newblock \emph{IEEE Transactions on Automatic Control}, 57(5), 1332--1338.

\bibitem[{Grip et~al.(2015)Grip, Fossen, Johansen, and Saberi}]{GripFJS2015AUT}
Grip, H.F., Fossen, T.I., Johansen, T.A., and Saberi, A. (2015).
\newblock Globally exponentially stable attitude and gyro bias estimation with
  application to {GNSS/INS} integration.
\newblock \emph{Automatica}, 51, 158--166.

\bibitem[{Higham(2008)}]{Higham2008book}
Higham, N.J. (2008).
\newblock \emph{Functions of Matrices: Theory and Computation}.
\newblock SIAM.

\bibitem[{Khalil(2002)}]{khalil2002book}
Khalil, H. (2002).
\newblock \emph{Nonlinear Systems}.
\newblock Prentice Hall.

\bibitem[{Mahony et~al.(2008)Mahony, Hamel, and Pflimlin}]{MahonHP2008TAC}
Mahony, R., Hamel, T., and Pflimlin, J.M. (2008).
\newblock Nonlinear complementary filters on the special orthogonal group.
\newblock \emph{IEEE Transactions on Automatic Control}, 53(5), 1203--1218.

\bibitem[{Martin and Sala{\"u}n(2007)}]{MartiS2007CDC}
Martin, P. and Sala{\"u}n, E. (2007).
\newblock Invariant observers for attitude and heading estimation from low-cost
  inertial and magnetic sensors.
\newblock In \emph{IEEE Conference on Decision and Control}, 1039--1045.

\bibitem[{Martin and Sala{\"u}n(2010{\natexlab{a}})}]{MartiS2010CEP}
Martin, P. and Sala{\"u}n, E. (2010{\natexlab{a}}).
\newblock {Design and implementation of a low-cost observer-based attitude and
  heading reference system}.
\newblock \emph{Control Engineering Practice}, 18(7), 712--722.

\bibitem[{Martin and Sala{\"u}n(2010{\natexlab{b}})}]{MartiS2010ICRA}
Martin, P. and Sala{\"u}n, E. (2010{\natexlab{b}}).
\newblock The true role of accelerometer feedback in quadrotor control.
\newblock In \emph{IEEE International Conference on Robotics and Automation},
  1623--1629.

\bibitem[{{Martin} and {Sarras}(2016)}]{MartiS2016arXiv}
{Martin}, P. and {Sarras}, I. (2016).
\newblock {A simple global observer for attitude and gyro biases}.
\newblock \emph{ArXiv e-prints}, arXiv:1604.03714v1 [math.OC].

\bibitem[{Tayebi et~al.(2013)Tayebi, Roberts, and Benallegue}]{TayebRB2013TAC}
Tayebi, A., Roberts, A., and Benallegue, A. (2013).
\newblock Inertial vector measurements based velocity-free attitude
  stabilization.
\newblock \emph{IEEE Transactions on Automatic Control}, 58(11), 2893--2898.

\bibitem[{Vasconcelos et~al.(2008)Vasconcelos, Silvestre, and
  Oliveira}]{VascoSO2008IFAC}
Vasconcelos, J., Silvestre, C., and Oliveira, P. (2008).
\newblock A nonlinear observer for rigid body attitude estimation using vector
  observations.
\newblock \emph{{IFAC} Proceedings Volumes}, 41(2), 8599 -- 8604.

\bibitem[{{Zamani} et~al.(2015){Zamani}, {Trumpf}, and
  {Mahony}}]{ZamanTM2015arXiv}
{Zamani}, M., {Trumpf}, J., and {Mahony}, R. (2015).
\newblock {Nonlinear Attitude Filtering: A Comparison Study}.
\newblock \emph{ArXiv e-prints}, arXiv:1502.03990 [cs.SY].

\end{thebibliography}
\end{document}